\documentclass[12pt]{amsart}

\AtBeginDocument{%
   \def\MR#1{}
}

\usepackage{amsmath,amsthm,amsfonts,amssymb,amscd,euscript,enumerate}
\usepackage[left=1in, right=1in, top=1.1in, bottom=1.1in]{geometry}
\usepackage{tikz}
\usepackage{tikz-cd}
\usetikzlibrary{fit}
\usetikzlibrary{trees}
\usepackage{forest}
\forestset{
  tree/.style={
    for tree={
      edge={line width=0.6pt},
      draw,
      line width=0.6pt,
      circle,fill=black,
      inner sep=0pt, minimum size=4pt,
      if n children=0{rectangle,fill=white}{}
    }
  }
}
\usepackage[all,arc]{xy}
\usepackage{enumerate}
\usepackage{mathrsfs} 
\usepackage{datetime}
\usepackage[colorlinks, citecolor = darkspringgreen, linkcolor = melon]{hyperref}
\usepackage{cleveref}
\usepackage{graphicx}
\usepackage{color}
\usepackage{comment}
\usepackage{float}
\input{xy}
\xyoption{all}
\usetikzlibrary{shapes}
\usetikzlibrary{tikzmark,calc}
\usepackage{mathtools}
\usepackage{listofitems}
\usepackage{wrapfig}
\usepackage{caption}
\usepackage{sidecap}
\usepackage{stmaryrd}
\usepackage{url}
\usetikzlibrary{decorations.pathreplacing, arrows.meta, calc}
\usepackage{enumitem}
\usetikzlibrary{patterns.meta}
\usepackage{pdfpages}

\definecolor{melon}{rgb}{0.4, 0.2, 1}
\definecolor{maroon}{rgb}{0.5,0,0}
\definecolor{violet}{rgb}{0.6, 0.3, 0.7}
\definecolor{darkviolet}{rgb}{0.35, 0, 0.6}
\definecolor{darkspringgreen}{rgb}{0.14, 0.6, 0.2}
\definecolor{amethyst}{rgb}{0.6, 0.4, 0.8}

\calclayout

\newcommand{\brect}[2]{%
  \begin{scope}[shift={(#1,#2)}]
    \fill[cyan] (0.2,0.2)--(0.2,2.8)--(0.8,2.8)--(0.8,0.2)--cycle;
  \end{scope}%
}

\newcommand{\rrect}[2]{%
  \begin{scope}[shift={(#1,#2)}]
    \fill[red] (-0.2,-0.2)--(-2.8,-0.2)--(-2.8,-0.8)--(-0.2,-0.8)--cycle;
  \end{scope}%
}

\setboolean{@twoside}{false}

\numberwithin{equation}{section}

\newtheorem{thm}{Theorem}[section]
\newtheorem{cor}[thm]{Corollary}
\newtheorem{prop}[thm]{Proposition}

\theoremstyle{definition}
\newtheorem{defn}[thm]{Definition}

\newtheorem{exmp}[thm]{Example}

\newtheorem*{conjs*}{Conjectures}

\theoremstyle{remark}
\newtheorem{rem}[thm]{Remark}

\crefname{thm}{theorem}{theorems}
\Crefname{thm}{Theorem}{Theorems}

\makeatletter
\let\c@equation\c@thm
\makeatother
\numberwithin{equation}{section}

\makeatletter
\newcommand*\bigcdot{\mathpalette\bigcdot@{.5}}
\newcommand*\bigcdot@[2]{\mathbin{\vcenter{\hbox{\scalebox{#2}{$\m@th#1\bullet$}}}}}
\makeatother

\makeatletter
\def\subsection{\@startsection{subsection}{3}%
  \z@{.5\linespacing\@plus.7\linespacing}{.1\linespacing}%
  {\bfseries}}
\makeatother

\newcommand{\Z}{\mathbb{Z}}
\newcommand{\N}{\mathbb{N}}

\newcommand{\calP}{\mathcal{P}}

\newcommand{\frakD}{\mathfrak{D}}

\renewcommand{\phi}{\varphi}

\DeclareMathOperator{\sh}{sh}

\DeclareMathOperator{\aux}{aux}
\DeclareMathOperator{\point}{point}

\newcommand{\east}{\mathbf{E}}
\newcommand{\north}{\mathbf{N}}
\newcommand{\edge}{e} 

\newcommand{\exc}[1]{\chi_{#1}} 
\newcommand{\LTS}[2]{\mathbf{LT}_{#2}^{#1}} 
\newcommand{\TaIn}[2]{\mathcal{I}_{#2}^{#1}} 
\newcommand{\TaLa}[2]{\mathcal{T}_{#2}^{#1}} 

\newcommand{\TGS}[2]{\mathbf{TG}_{#2}^{#1}} 
\newcommand{\TG}[3]{\mathbf{TG}_{#2,#3}^{#1}} 
\newcommand{\PLTS}[2]{\rlap{$\widehat{\mathbf{LT}}$}\phantom{\mathbf{LT}}_{#2}^{#1}} 

\title{{Bijecting \lowercase{$m$}-Tamari intervals, trees, and tight gradings}}

\author{Amanda Burcroff}
\address{Amanda Burcroff, Department of Mathematics, MIT, Cambridge, MA 02139, USA}
\email{\href{mailto:amandabu@mit.edu}{amandabu@mit.edu}}

\author{Kyungyong Lee}
\address{Kyungyong Lee, Department of Mathematics, University of Alabama, Tuscaloosa, AL 35487, USA and Korea
Institute for Advanced Study, Seoul 02455, Republic of Korea}
\email{\href{mailto:kyungyong.lee@ua.edu}{kyungyong.lee@ua.edu} \& \href{mailto:klee1@kias.re.kr}{klee1@kias.re.kr}}

\author{Lang Mou}
\address{Lang Mou, Department of Mathematics, University of California Davis, Davis, CA 95616, USA}
\email{\href{mailto:lmou.math@gmail.com}{lmou.math@gmail.com}}

\date{\today}

\begin{document}

\begin{abstract}
    Tight gradings, which are certain integer gradings on maximal Dyck paths, provide a combinatorial framework for computing all coefficients of wall-functions in rank-$2$ cluster scattering diagrams.  We construct a bijection between a special subset of tight gradings and the intervals of the $m$-Tamari lattice, providing a new combinatorial proof of the central wall-function formula proved by Reineke. Inspired by the structure of these gradings, we introduce a class of plane trees called \emph{$m$-leafy trees} that give a new model for $m$-Tamari intervals.  The defining condition of an $m$-leafy tree is simple: in every full subtree, the number of leaves is greater than $m$ times the number of internal vertices.  The leafy tree perspective gives rise to a quadratic functional equation with a catalytic parameter, providing an alternative way to enumerate $m$-Tamari intervals. 
\end{abstract}

\keywords{plane tree, $m$-Tamari lattice, wall-crossing, tight grading, cluster algebra}


\maketitle

\vspace{-1em}

\section{Introduction}

The Tamari lattice $\mathcal{T}_{n}$, introduced by Tamari in 1951, is the partial order on binary trees with $n$ internal vertices generated by rotations. Its Hasse diagram is an orientation of the edge graph of an associahedron, linking Catalan combinatorics with the geometry of associativity. Tamari intervals, namely comparable pairs in $\mathcal{T}_n$, gained sustained interest after Chapoton showed in 2006 \cite{Chapoton05} that they are counted by $\frac{2}{n(n+1)}\binom{4n+1}{n-1}$. This simple formula inspired the development of several combinatorial models, for example \cite{BB,ChatelPons2015,FFN}.

Bergeron observed that Chapoton’s formula coincides with Haiman’s conjectural dimension of the alternating component of the $\mathbb{S}_n$-module of trivariate diagonal harmonics. Motivated by $m$-variant versions of these spaces, he introduced the $m$-Tamari lattice $\mathcal{T}_{n}^{m}$ \cite{BPR} and conjectured that its intervals are counted by
\begin{equation}\label{eq:number of m-Tamari int intro}
    \frac{m+1}{n(mn+1)}\binom{(m+1)^2n+m}{n-1}.
\end{equation}

Bousquet-M\'elou, Fusy, and Pr\'eville-Ratelle subsequently proved this formula counting the $m$-Tamari intervals \cite{BMFPR}.

More recently Chapoton observed that the numbers \eqref{eq:number of m-Tamari int intro} also coincide with the Euler characteristics of a family of Kronecker moduli spaces \cite{ACMR}, which were computed by Weist \cite{Weist}. These Kronecker moduli are known to be closely related to scattering diagrams in rank $2$ \cite{RW, GP, BLMMR}. 

Our approach to $m$-Tamari intervals arises from the computation of cluster scattering diagrams \cite{GHKK}. In prior work \cite{BLMshort}, the authors introduced \emph{tight gradings}, which provide a uniform combinatorial formula for wall-function coefficients in any rank-2 cluster scattering diagram. These gradings on Dyck paths are built on Lee--Li--Zelevinsky’s compatible pairs \cite{LLZ}, and play a crucial role in extending Laurent positivity to generalized cluster algebras \cite{BLM}. While seeking a combinatorial enumeration of tight gradings associated with the central wall, we noticed that the number of $0$-framed tight gradings is precisely that of $m$-Tamari intervals. This equinumerosity motivated us to look for an explicit combinatorial connection.

To establish this connection, we introduce \emph{$m$-leafy trees} as a novel plane tree model for $m$-Tamari intervals. Their defining condition is surprisingly simple: in every full rooted subtree, the number of leaves is strictly larger than $m$ times the number of internal vertices. These trees serve as intermediate objects bridging $m$-Tamari intervals and tight gradings.

The main aim of this extended abstract is to present the bijections between these three sets of combinatorial objects:
\begin{enumerate}[label=(\Roman*)]
    \item intervals in the $m$-Tamari lattice of order $n$,
    \item \label{item: m-leafy} $m$-leafy trees on $n$ internal vertices and $mn+1$ leaves, and
    \item $0$-framed tight gradings contributing to the $n^{\text{th}}$ central wall coefficient of the cluster scattering diagram $\frakD_{m+2,m+2}$.
\end{enumerate}

\begin{figure}[b]
\centering
\includegraphics[width=\textwidth]{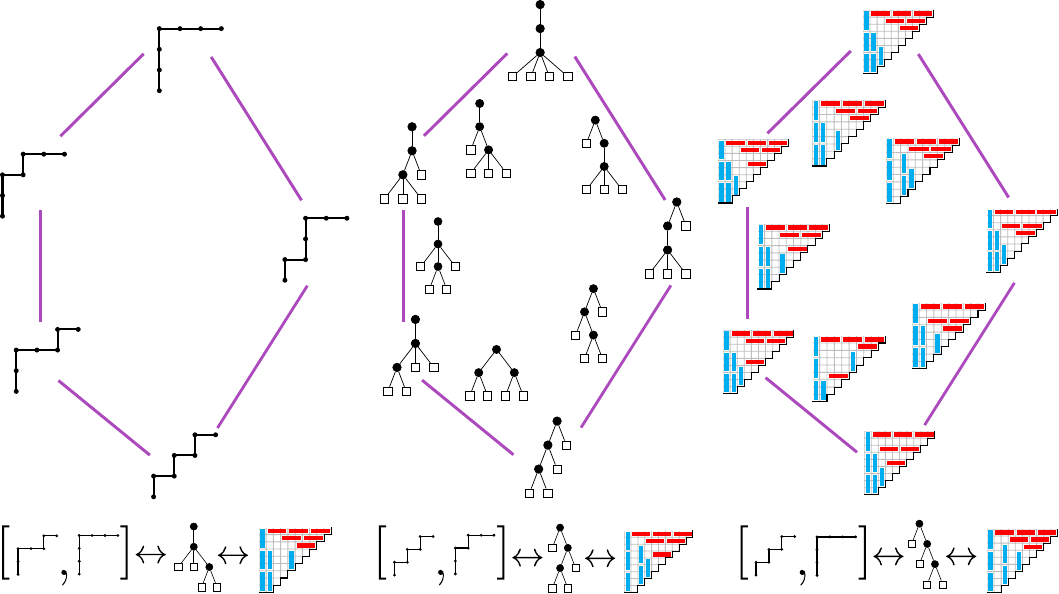}
\caption{For the Hasse diagram of the Tamari lattice $\mathcal{T}_3 = \mathcal{T}_{3}^{1}$ (left), we depict the leafy trees (center) and tight gradings (right) corresponding to the vertices and edges.  The intervals of size greater than $2$ are described below.}
\label{fig:introduction-bijections}
\end{figure}

The corresponding $m$-Tamari interval in (I) can be read directly from the right-to-left preorder of a leafy tree in (II): the maximum element records whether each vertex is internal or a leaf, while the minimum element records the arities of the internal vertices. Moreover, we extend the bijection between (II) and (III) to a bijection between $m$-leafy trees with a distinguished vertex and the set of all tight gradings contributing to the $n^{\text{th}}$ central wall coefficient of $\frakD_{m+2,m+2}$.

 The bijection to $m$-leafy trees yields several nice consequences for $m$-Tamari intervals.  The recursive structure of leafy trees makes them easy to generate and enumerate, and they yield a comparability test for $m$-Tamari elements as a positivity condition on plane trees. We find that a generating series for $m$-leafy trees, refined by two tree statistics, satisfies the functional equation \eqref{eq:functional equation 1} obtained in \cite{BMFPR} for $m$-Tamari intervals. We also derive a new quadratic functional equation \eqref{eq:functional equation 2} with one catalytic parameter using the tree structure. Solving this equation gives an alternative simpler derivation of the formula \eqref{eq:number of m-Tamari int intro}. The obvious left-right symmetry on $m$-leafy trees induces an involution on $m$-Tamari intervals which was not previously identified, to the best of our knowledge.

The connection with tight gradings gives a second application. Although tight gradings are proven to compute wall coefficients \cite{BLM}, we had not been able to combinatorially enumerate them, even for the central wall, until linking them with leafy trees. So the tree perspective provides a new proof of the central wall formula for $\mathfrak{D}_{m+2,m+2}$, conjectured by Kontsevich in a closely related setting of tropical vertex \cite{GP} and proved by Reineke using quiver representations \cite{Rei11, BLMMR}.

\section{\texorpdfstring{A plane-tree model for $m$-Tamari intervals}{A tree model for m-Tamari intervals}}

In this section, we review the $m$-Tamari lattice and its intervals, introduce $m$-leafy trees, and describe the bijection $\Phi$ from $m$-leafy trees to $m$-Tamari intervals.

\subsection{\texorpdfstring{Intervals in the $m$-Tamari lattice}{Intervals in the m-Tamari lattice}}

Dyck paths of size $n$ form the Tamari lattice $\mathcal{T}_n$. An \emph{$m$-ballot path} of size $n$ is a word with $n$ letters $N$ and $mn$ letters $E$ such that every prefix $W$ satisfies $m|W|_N - |W|_E \geq 0$, where $|W|_L$ is the number of letter $L$ in $W$. Equivalently, it is a north-east lattice path from $(0,0)$ to $(mn, n)$ that stays non-strictly above $x=my$. Thus $m$-ballot paths generalize Dyck paths where $m=1$.

\begin{defn}
    For $P$ and $Q$ two $m$-ballot paths of size $n$, define $P\prec Q$ if there exists in $P$ an east step $a$, followed by a north step $b$, such that $Q$ is obtained from $P$ by swapping $a$ and $S$, where $S$ is the shortest factor of $P$ that begins with $b$ and is a (translated) $m$-ballot path.
\end{defn}

The transitive closure $\leq$ of $\prec$ defines a lattice on $m$-ballot paths of size $n$, which we call the \emph{$m$-Tamari lattice} $\TaLa{m}{n}$, generalizing $\TaLa{1}{n} = \mathcal{T}_n$. Any pair $P, Q\in \TaLa{m}{n}$ with $P\leq Q$ is an \emph{interval} and we denote it by $[P, Q]$. Write $\TaIn{m}{n}$ for the set of all intervals.

We refer the reader to \cite{BMFPR} for a more comprehensive exposition on $m$-Tamari lattices and their intervals. We note that there have been combinatorial models for $m$-Tamari intervals, for example \cite{ChatelPons2015,Pons19}, and it would be interesting to understand the relation with our tree model introduced below.

\subsection{\texorpdfstring{$m$-leafy trees}{m-leafy trees}}

All trees considered here are finite plane trees; in particular, there is a designated root and the children of each \emph{internal vertex} (i.e. non-leaf vertex) are ordered from left to right. In our illustrations, the root is always placed at the top.

A \emph{subtree} of a tree $T$ is a plane tree rooted at a vertex $v$ of $T$ spanned by all descendants of $v$, and is denoted by $T_v$.

Let $\ell(T)$ denote the number of leaves of $T$ and $n(T)$ the number of internal vertices.

Fix a positive integer $m$. The following notion of \emph{$m$-leafy trees} is central to our discussion but remarkably simple.

\begin{defn}[$m$-leafy trees]
    A tree $T$ is called \emph{$m$-leafy} if for every vertex $v$,
    \[
        \exc{m}(T_v) \coloneqq \ell({T_v}) - mn({T_v}) > 0.
    \]
\end{defn}

Despite this simple characterization, the class of $m$-leafy trees seems not to have drawn much attention in combinatorics. However, some special cases are familiar families. For example, $m$-leafy trees with $\exc{m}(T_v)=1$ for every $v$ are precisely $(m+1)$-ary trees.

\begin{figure}[ht]
    \centering
    \begin{forest}
      baseline,
      tree,
      /tikz/every label/.append style={font=\footnotesize},
      for tree={
        l=7mm, l sep=0mm, s sep=3mm,
      }
      [,label=right:$1$
        [, label=left:$1$]
        [, label=right:$2$
          [, label=left:$3$
            [, label=below:$1$] [, label=below:$1$] [, label=below:$1$] [, label=below:$1$] [, label=below:$1$] 
          ]
          [, label=right:$1$]
        ]
      ]
    \end{forest}
    \caption{A $2$-leafy tree with $\exc{2}(T_v)$ next to each vertex $v$.}
    \label{fig:2-leafy with 2-excess}
\end{figure}
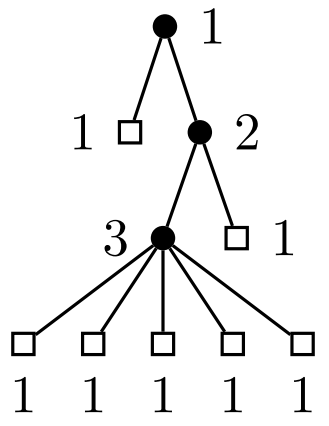

For $n\geq 1$, write $\LTS{m}{n} \coloneqq \left\{T\colon \text{$T$ is $m$-leafy with $n$ internal vertices and $\exc{m}(T)=1$}\right\}$.  Equivalently, $\LTS{m}{n}$ is the set of $m$-leafy trees with $n$ internal vertices and $mn+1$ leaves.

\subsection{\texorpdfstring{Bijecting $m$-leafy trees and $m$-Tamari intervals}{Bijecting m-leafy trees and m-Tamari intervals}}

The promised bijection $\Phi \colon \LTS{m}{n} \to \TaIn{m}{n}$ is easy to describe.

Let $S$ be a tree. The \emph{right-to-left preorder} of $S$ is defined to visit the root first, then recursively visit the children from right to left. List among all vertices the internal vertices in this order $v_1, v_2, \dots, v_n$.
Let $k_i$ be the \emph{arity} of $v_i$, the number of children of $v_i$.

We associated the \emph{upper} and \emph{lower} words to $S$:
\begin{align}
    U_S &\coloneqq
    \begin{array}[t]{@{}l@{}}
    \text{the word of all vertices in the right-to-left preorder,}\\
    \text{with $N$ for internal vertices and $E$ for leaves, after deleting the final $E$;} 
    \end{array}\label{eq:upper path}\\
    L_S &\coloneqq NE^{k_1-1}NE^{k_2-1}\cdots NE^{k_n-1}. \label{eq:lower path}
\end{align}

\begin{figure}[ht]
\centering
\begin{tabular}{@{}c@{\hspace{12mm}}l@{}}
\begin{forest}
  baseline,
  tree,
  for tree={
    l=6mm,l sep=0pt,s sep=2.5mm,
  }
  [
    []
    [
      [ [] [] [] [] [] ]
      []
    ]
  ]
\end{forest}
&
\raisebox{-12mm}
{$
  \begin{aligned}
    U_S &= {NNENEEEEE}\\[4mm]
    L_S &= {NENENEEEE}
  \end{aligned}
$}
\\ \noalign{\vskip 6pt}
$S$ &
\end{tabular}
\caption{A $2$-leafy tree in $\LTS{2}{3}$ and the associated words.}
\label{fig:eg 2-leafy tree and paths}
\end{figure}

\begin{thm}\label{thm:general m bijection}
    If $S\in \LTS{m}{n}$, then $L_S$ and $U_S$ are $m$-ballot paths with $L_S\leq U_S$. Moreover, the map below is a bijection:
    \begin{equation}\label{eq:m-tree to m-interval bijection}
        \Phi \colon \LTS{m}{n} \to \TaIn{m}{n},\qquad \Phi(S) \coloneqq [L_S, U_S].
    \end{equation}
\end{thm}

Next we discuss some consequences of the above tree model of $m$-Tamari intervals.

A plane tree is called a \emph{caterpillar} if its internal vertices form a chain. A caterpillar is called a \emph{left} (resp. \emph{right}) caterpillar if any internal vertex is the leftmost (resp. rightmost) child of its parent.

The bijection $\Phi$ gives rise to three models for the $m$-Tamari lattice:
\begin{itemize}
    \item Left $m$-leafy caterpillars are length-zero intervals, namely $m$-ballot paths;
    \item Right $m$-leafy caterpillars are intervals whose upper path is maximal;
    \item ($m$+1)-ary trees are intervals whose lower path is minimal.
\end{itemize}

\begin{figure}[htbp]
\centering
\begin{forest}
  baseline,
  tree,
  for tree={
    l=6mm,l sep=0pt,s sep=2.5mm,
  }
  [
    [
      [
        [ [] [] [] ]
        [] [] [] []
      ]
      []
    ]
    []
  ]
\end{forest}
\hspace{10mm}
\begin{forest}
  baseline,
  tree,
  for tree={
    l=6mm,l sep=0pt,s sep=2.5mm,
  }
  [
    []
    [
      []
      [
        [] [] [] []
        [ [] [] [] ]
      ]
    ]
  ]
\end{forest}
\hspace{10mm}
\begin{forest}
  baseline,
  tree,
  for tree={
    l=8mm,l sep=0pt,s sep=2.5mm,
  }
  [
    [ [] [] [] ]
    [
      []
      [ [] [] [] ]
      []
    ]
    []
  ]
\end{forest}
\caption{A left caterpillar, a right caterpillar, and a ternary tree corresponding to the same $2$-ballot path ${NENENEEEENEE}$.}
\label{fig:three models m-Tamari}
\end{figure}
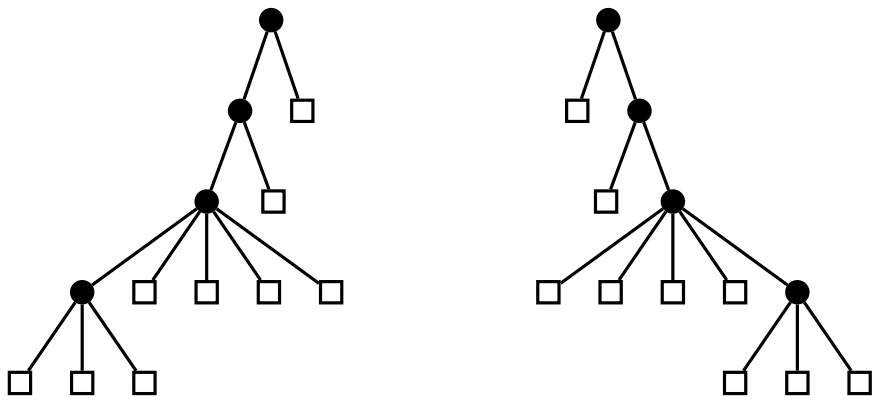

The third model is a well-known $(m+1)$-ary tree description of the $m$-Tamari lattice. The cover relation is realized as certain tree rotation, which was already pointed out in \cite[Section 4]{BMFPR}.

\begin{figure}
\centering
\begin{forest}
  baseline,
  tree,
  for tree={
    l=0pt, l sep=8mm,s sep=3mm,
  }
  [,tikz={
    \node[below=5pt] at (current bounding box.south) {$S$};
  }
    []
    [,label=right:$v$
      [ [] [] [] [] [] ]
      []
    ]
  ]
\end{forest}
\quad
\raisebox{-16mm}{$\overset{\Phi}{\longmapsto}$}
\quad
\begin{forest}
  baseline,
  tree,
  for tree={
    l=0pt, l sep=8mm,s sep=7mm,
  }
  [
    []
    [,label=right:$v$
      [,tikz={
        \node[draw,dashed,rectangle,rounded corners=5pt,
              inner sep=3pt,label distance=2pt,
              fit=(),label=below:$T_2$] {};
      }]
      [,s sep=2mm,tikz={
        \node[draw,dashed,rectangle,rounded corners=5pt,
              inner sep=3pt,label distance=2pt,
              fit=() (!1) (!2) (!3),label=below:$T_1$] {};
      }
        [] [] []
      ]
      [,tikz={
        \node[draw,dashed,rectangle,rounded corners=5pt,
              inner sep=3pt,label distance=2pt,
              fit=(),label=below:$T_0$] {};
      }]
    ]
    [,label=right:$u$]
  ]
\end{forest}
\quad
\raisebox{-16mm}{$\prec$}
\quad
\begin{forest}
  baseline,
  tree,
  for tree={
    l=0mm, l sep=8mm,s sep=7mm,
  }
  [
    []
    [,tikz={
      \node[draw,dashed,rectangle,rounded corners=5pt,
            inner sep=3pt,label distance=2pt,
            fit=(),label=below:$T_2$] {};
    }]
    [,label=right:$v$
      [,label=below:$u$]
      [,s sep=2mm,tikz={
        \node[draw,dashed,rectangle,rounded corners=5pt,
              inner sep=3pt,label distance=2pt,
              fit=() (!1) (!2) (!3),label=below:$T_1$] {};
      }
        [] [] []
      ]
      [,tikz={
        \node[draw,dashed,rectangle,rounded corners=5pt,
              inner sep=3pt,label distance=2pt,
              fit=(),label=below:$T_0$] {};
      }]
    ]
  ]
\end{forest}
\caption{A $2$-leafy tree $S$ that is a cover relation.}
\label{fig:2-leafy cover relation}
\end{figure}
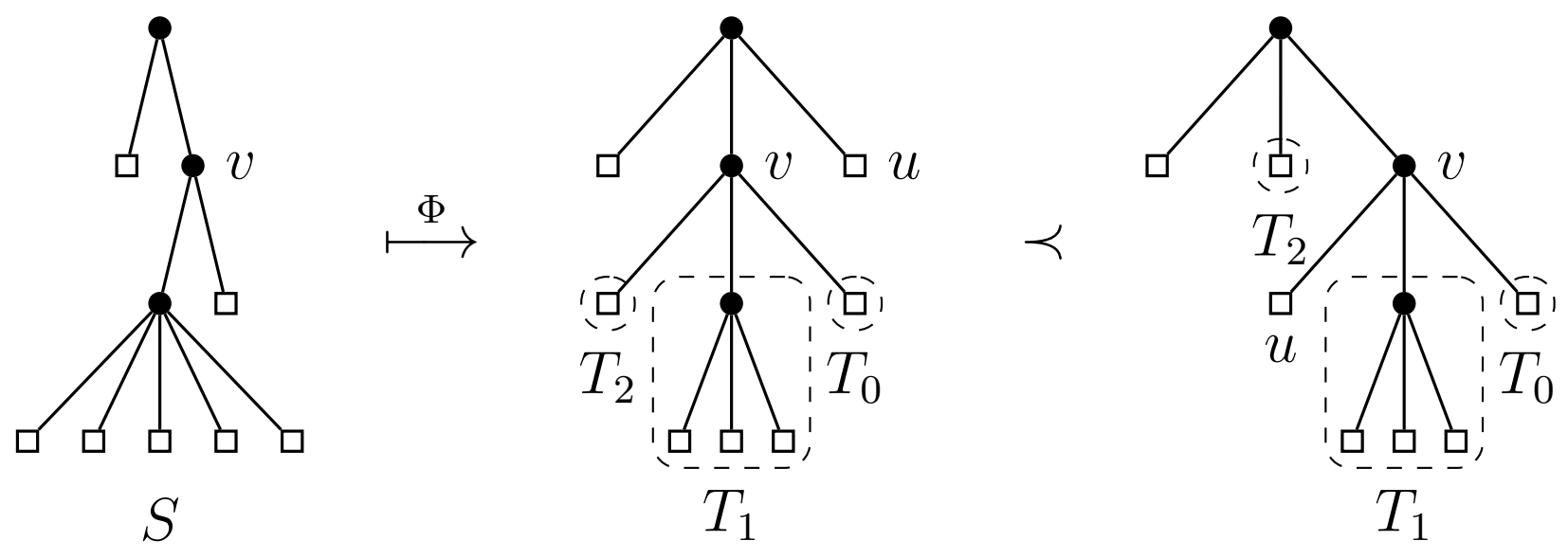

We may characterize the $m$-leafy trees that correspond to cover relations in $\TaLa{m}{n}$. See \Cref{fig:2-leafy cover relation} for an example.

\begin{prop}
    For $S\in \LTS{m}{n}$, the interval $\Phi(S) = [L_S, U_S] \in \TaIn{m}{n}$ is a cover relation if and only if the two conditions hold:
    \begin{itemize}[font=\upshape]
        \item $S$ has precisely one internal vertex $v$ that is the second leftmost child of its parent, and any other non-root internal vertex is the leftmost child of its parent.
        \item $\exc{m}(S_v) < \exc{m}(S_v')$ for any internal descendant $v'$ of $v$.
    \end{itemize}
\end{prop}

The left-right reflection of plane trees defines an involution $\tau\colon \LTS{m}{n} \to \LTS{m}{n}$.
In particular, it sends a left caterpillar to its mirrored right caterpillar, and vice versa. The set of $(m+1)$-ary trees is closed under $\tau$.

The involution $\tau$ gives rise to the involution on $m$-Tamari intervals
\[
    \gamma_m\coloneqq \Phi \circ \tau \circ (\Phi)^{-1}\colon \TaIn{m}{n} \to \TaIn{m}{n},
\]
which exchanges length-zero intervals and upper intervals. We remark that $\gamma_m$ is different from the involution $\beta_m$ constructed by Pons since $\beta_m$ sends length-zero intervals to length-zero ones \cite[(5.36)]{Pons19}.

The last page after the bibliography illustrates \Cref{thm:general m bijection} for $\Phi\colon \LTS{2}{3} \to \TaIn{2}{3}$, where $2$-ballot paths are represented by left $2$-leafy caterpillars.

\section{\texorpdfstring{Enumeration of $m$-leafy trees}{Enumeration of m-leafy trees}}

The bijection $\Phi\colon \LTS{m}{n} \to \TaIn{m}{n}$ already implies the closed-form expression
\begin{equation}\label{eq:number of leafy trees}
    |\LTS{m}{n}| = \frac{m+1}{n(mn+1)} \binom{(m+1)^2n+m}{n-1},
\end{equation}
as $|\TaIn{m}{n}|$ equals the same expression \eqref{eq:number of m-Tamari int intro}.

Trees are among typical combinatorial objects that can be constructed recursively based on disjoint union and concatenation. This feature makes it possible to write down functional equations of their generating series. We present two functional equations involving generating series of $m$-leafy trees. The first already appeared in the enumeration of $m$-Tamari intervals \cite{BMFPR}. The second appears to be new and can be solved following general methods developed in \cite{BMJ}; see also \cite{BM06}.

\subsection{First functional equation}\label{subsec:first func eq}

For a tree $T$, let $l(T)$ (resp. $r(T)$) be the number of internal vertices in the longest chain obtained by taking only the leftmost (resp. rightmost) child starting from the root. For a single leaf, $l = r = 0$. Consider the generating function
\begin{equation}\label{eq:gen function left right lengths}
    F^{(m)}(t;x,y) \coloneqq \sum_{n\geq 0}\sum_{T\in \LTS{m}{n}} t^n x^{l(T)+1}y^{r(T)}.
\end{equation}

\begin{prop}\label{prop:functional equation 1}
    Writing $F(x,y)\equiv F^{(m)}(t;x,y)$, we have
    \begin{equation}\label{eq:functional equation 1}
        F(x,y) = x + xyt\left(F(x,1)\cdot \Delta\right)^m\left(F(x,y)\right),
    \end{equation}
    where $\Delta$ is the divided difference operator $\Delta S(x) = \frac{S(x)-S(1)}{x-1}$.
\end{prop}

The equation \eqref{eq:functional equation 1} is precisely \cite[(3)]{BMFPR}, although interpreted in terms of $m$-Tamari intervals using the initial rise of upper path and the contact of lower path. We remark that via the bijection $\Phi$, the initial rise translates to our $r$ whereas the contact does not translate to our $l$. In fact, the quantity $l(T)$ is the number of north steps shared by the $m$-ballot paths $L_T$ and $U_T$, while $r(T)$ is the number of initial north steps of $U_T$.

The tree interpretation makes the following symmetry (described as a \emph{combinatorial mystery} for rise-contact in \cite{BMFPR}) evident, since the mirror map $\tau$ exchanges the statistics $l$ and $r$:
\begin{equation}\label{eq:left/right symmetry}
    F^{(m)}(t; x, y)/{x} = F^{(m)}(t; y, x)/{y}.
\end{equation}

\subsection{Second functional equation}\label{subsec:second func eq}

While the solution to \eqref{eq:functional equation 1} in \cite{BMFPR} is rather complex, we present another functional equation with a \emph{catalytic parameter} $u$ whose solution is simpler. Define the generating series
\begin{equation}\label{eq:gen fun forest with excess}
    G^{(m)}(t, u) = 1 + \sum_{T} t^{n(T)} u^{\chi_m(T)} \in \mathbb{Q}[[t, u]],
\end{equation}
where the sum is over all forests $T$ whose component trees are $m$-leafy, and $n(T)$ and $\chi_m(T)$ are extended from trees to forests additively.

Write $G = G^{(m)}(t,u)$ and $G_j \equiv G_j^{(m)}(t) \coloneqq [u^j]G^{(m)}(t,u)$ for $j\geq 0$, where $G_0=1$. For example, $G_j(0)=1$, which counts the forest of isolated $j$ leaves; $G_1$ is the generating series of the family $\{\LTS{m}{n}\}_{n\geq 0}$.

\begin{prop}
    The series $G$ and $G_j$, $j=1,\dots, m$ satisfy the functional equation
    \begin{equation}\label{eq:functional equation 2}
        tG^2 - \left( u^m(1-u) + t \sum_{j=0}^{m} G_j u^j \right) G + u^m = 0.
    \end{equation}
\end{prop}

Note that \eqref{eq:functional equation 2} determines $G$ and $G_1,\dots, G_m$ all at once. Solving for $G_1$ following the general methods in \cite{BMJ} leads to the closed formula \eqref{eq:number of leafy trees} of $|\LTS{m}{n}|$.

\section{Bijection to central tight gradings}

While leafy trees are a new combinatorial object to our knowledge, they arose naturally in the study of tight gradings recently introduced by the authors to investigate rank-2 scattering diagrams \cite{BLMshort}.
In this section, we construct an explicit bijection between leafy trees and certain tight gradings.  We begin by providing preliminaries on the relevant tight gradings, which are those corresponding to the central wall-coefficients of rank-$2$ skew-symmetric cluster scattering diagrams. Cluster algebras and their associated cluster scattering diagrams have an incredibly rich theory with numerous connections across mathematics and physics which we do not describe here; we instead refer the reader to \cite{FWZ, GHKK, NakTextbook} for further details.  

\subsection{Central tight gradings}\label{subsec:tight grading}

We now give an explicit description of the tight gradings corresponding to the $n^\text{th}$ central wall coefficient of the cluster scattering diagram $\frakD_{m+2,m+2}$.  Tight gradings are defined much more generally in \cite{BLM}, and a tiling description is given in \cite[Section 5]{BLMMR}.

For $m,n \in \N$, let $\calP_{m,n}$ be the $NE$-path $E(EN)^{(m+2)n}$.  This path is sometimes referred to as the \emph{maximal Dyck path} inside the rectangle of length $(m+2)n+1$ and height $(m+2)n$.  Given a subset $C$ of steps (also known as \emph{edges}) in an $NE$-path $\calP$, we denote the set of east (resp. north) steps by $C_\east$ (resp. $C_\north$) and the total number of steps by $|C|$.

A \emph{grading} on an $NE$-path is an assignment of a nonnegative integer to each edge of the path.   Given a grading $\omega: \mathcal{P} \to \Z_{\ge 0}$ and a set of edges $C \subset \mathcal{P}$, we set $\omega(C) \coloneqq \sum_{\edge \in C}\omega(\edge)$.   For edges $u \in \calP_\east$ and $v \in \calP_\north$, let $\overrightarrow{uv}$ be the subpath proceeding east from the left endpoint of $u$ to the top endpoint of $v$, continuing cyclically around $\calP$ if necessary.

\begin{defn}\label{def: shadow}
Fix an $NE$-path $\calP$ and a grading $\omega \colon \calP \rightarrow \Z_{\geq 0}$.
For each edge $\edge$ in $\calP$, we define its \emph{shadow}, denoted by $\sh(\edge)$ or $\sh(\edge;\omega)$, as follows. 
\begin{itemize}
    \item If $\edge \in \calP_\east$, then its shadow is $\overrightarrow{\edge v}_\north$, where $v \in \calP_\north$ is chosen such that ${|\overrightarrow{\edge v}_\north| = \omega(\overrightarrow{\edge v}_\east)}$ and $\overrightarrow{\edge v}$ has minimal length under this condition.  If no such $v$ exists, let $\sh(\edge) = \calP_\north$.
    \item If $\edge \in \calP_\north$, then its shadow is $\overrightarrow{u\edge}_\east$, where $u \in \calP_\east$ is chosen such that $|\overrightarrow{u\edge}_\east| = \omega(\overrightarrow{u\edge}_\north)$ and $\overrightarrow{u\edge}$ has minimal length under this condition.  If no such $u$ exists, let $\sh(\edge) = \calP_\east$.
\end{itemize}

For $C \subset \calP$, let the \emph{shadow} of $C$ be $\sh(C)=\bigcup_{\edge \in C} \sh(\edge)$.  
\end{defn}

\begin{defn}
    A grading $\omega$ is \emph{compatible} if for any edges $u \in \calP_\east$ and $v \in \calP_\north$, either 
    $u \notin \sh(v)$ or $v \notin \sh(u)$.
\end{defn}

Compatible gradings (or compatible pairs) originated in the work of Lee-Li-Zelevinsky in their construction of the greedy basis \cite{LLZ}.

\begin{defn}
    A \emph{central tight grading of order $n$} is a compatible grading on $\calP_{m,n}$ such that there are $n$ horizontal and $n$ vertical edges having value $m+2$ in the grading, and all other edges have value $0$ in the grading.  We denote the set of all such tight gradings by $\TGS{m}{n}$.
\end{defn}

The following result is a special case of \cite[Theorem 2.3]{BLMshort}.

\begin{thm}
    The central wall-function in $\mathfrak{D}_{m+2,m+2}$ is given by
    \[
        1 + \sum_{n\geq 1} |\TGS{m}{n}| x^{(m+2)n}y^{(m+2)n}.
    \]
\end{thm}

Pictorially, we represent gradings by drawing a stack of height $|\sh(u)|$ consisting of $(m+2) \times 1$  blue rectangular tiles above each $u \in \calP_\east$, and a row of length $|\sh(v)|$ consisting of $1 \times (m+2)$ red rectangular tiles to the left of each $v \in \calP_\north$.  Given a path $\calP = UNE^kN$ where $U$ is some $NE$-word, we define the \emph{rotation operation} $\varrho$ which takes $\calP$ to be the $NE$-path $\varrho(\calP) \coloneqq E^kNUN$. A \emph{rotation} of $\calP$ is any path of the form $\varrho^r(\calP)$ for $r \geq 0$.  For the pictorial representation of a grading, it is sometimes useful to instead work on a rotation of the underlying $NE$-path.

\begin{figure}[ht]
\centering
\begin{tikzpicture}[scale=.35,
   labelled/.style={
    circle, draw, minimum size=6mm, inner sep=0pt
  },
  dot/.style={
    circle, fill=black, inner sep=0pt, minimum size=7pt
  },
  level distance=25mm,
  sibling distance=30mm,
  edge from parent/.style={draw}
]
        \begin{scope}
        \draw[step=1,color=gray] (0,0) grid (10,9);
        \end{scope}
        \draw[line width=1.5,color=black] (0,0)--(2,0)--(2,1)--(3,1)--(3,2)--(4,2)--(4,3)--(5,3)--(5,4)--(6,4)--(6,5)--(7,5)--(7,6)--(8,6)--(8,7)--(9,7)--(9,8)--(10,8)--(10,9);
        \brect{0}{0}
        \brect{0}{3}
        \brect{0}{6}
        \brect{1}{0}
        \brect{1}{3}
        \brect{3}{2}
        \rrect{10}{9}
        \rrect{7}{9}
        \rrect{4}{9}
        \rrect{9}{8}
        \rrect{6}{8}
        \rrect{7}{6}
    \begin{scope}[shift={(16,0)}]
    \begin{scope}
    \draw[step=1,color=gray] (0,0) grid (10,9);
        \end{scope}
        \draw[line width=1.5,color=black] (0,0)--(1,0)--(1,1)--(2,1)--(2,2)--(3,2)--(3,3)--(4,3)--(4,4)--(5,4)--(5,5)--(6,5)--(6,6)--(7,6)--(7,7)--(9,7)--(9,8)--(10,8)--(10,9);
        \brect{0}{0}
        \brect{0}{3}
        \brect{0}{6}
        \brect{1}{1}
        \brect{1}{4}
        \brect{3}{3}
        \rrect{10}{9}
        \rrect{7}{9}
        \rrect{4}{9}
        \rrect{9}{8}
        \rrect{6}{8}
        \rrect{7}{7}
    \end{scope}
    \end{tikzpicture} 
    \caption{A $0$-framed central tight grading (left) and a $7$-framed central tight grading (right) where $m = 1$ and $n = 3$.  Note that edges of the underlying $NE$-path that are directly below a blue tile or to the right of a red tile are precisely the edges with value $m+2 = 3$ in the tight grading, and all other edges take value $0$.} 
    \label{fig: tight grading examples}
\end{figure}
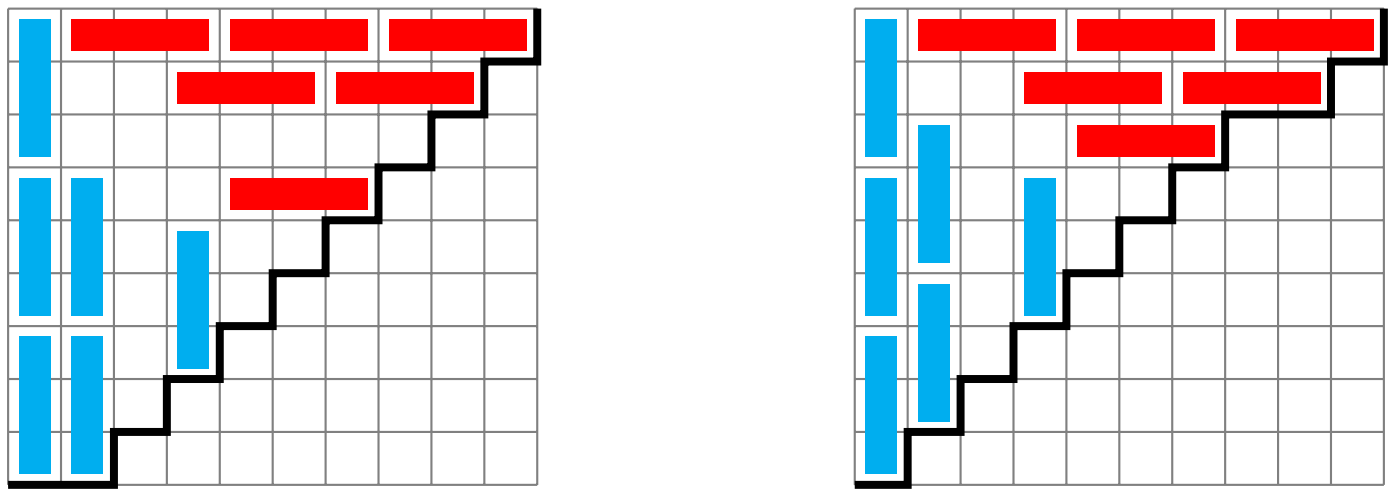

For each central tight grading $\omega$ on $\calP_n$, there is a unique choice of $f \in \{0,1,\dots,(m+2)n-1\}$ such that all tiles in the pictorial representation of the grading are contained in the rectangle bounding the rotated path $\varrho^{f}$ \cite[Remark 5.2]{BLMMR}.  In this case, we say that $\omega$ is \emph{$f$-framed}.

\begin{rem} 
The compatibility condition on tight gradings implies that in this pictorial representation, the red and blue tiles do not overlap when drawn on any rotation of the underlying Dyck path.
\end{rem}

\subsection{Bijection between tight gradings and  vertex-pointed leafy trees}

A tree is \emph{vertex-pointed} if it is equipped with a distinguished choice of vertex, denoted by $v_{\point}$.  Let $\PLTS{m}{n}$ denote the set of vertex-pointed $m$-leafy trees on $n$ internal vertices and $mn+1$ leaves.  Note that $\PLTS{m}{n} \cong \LTS{m}{n} \times [mn + n + 1]$. 

Fix $T \in \PLTS{m}{n}$.  For each internal vertex $v$, add an auxiliary vertex to be the rightmost child of $v$.  This yields a new tree $\aux(T)$, which is an $(m+1)$-leafy tree with $mn + n + 1$ leaves and $n$ internal vertices.  Note that the distinguished vertex $v_{\point}$ in $\aux(T)$ can be any vertex except the auxiliary vertices.

Label all vertices of $\aux(T)$ according to the left-to-right preorder, denoting the label of a vertex $v$ by $\lambda(v)$ and starting with label $0$ at the root.  We will work on the rotated NE path $\varrho^f(\calP_n)$, where $f \coloneqq \lambda(v_{\point})$.

For each horizontal (resp. vertical) edge of $\varrho^f(\calP_n)$, label it by $u_j$ (resp. $v_{j-1}$) where $j$ is the number of horizontal edges to its left in $\varrho^f(\calP_n)$.  Let $\omega_\east$ denote the set of horizontal edges with nonzero value in the grading $\omega$, and analogously for $\omega_\north$.

\begin{thm}\label{thm: leafy to tight bij}
There is an explicit bijection
\begin{equation}\label{eq:leafy to tight}
    \widehat \Gamma \colon \PLTS{m}{n} \to \TGS{m}{n}
\end{equation}
defined as follows.  Given a vertex-pointed leafy tree $(T,v_{\point}) \in \PLTS{m}{n}$ with $\lambda(v)$ denoting the label of each vertex $v \in \aux(T)$, the corresponding tight grading $\omega \coloneqq \widehat \Gamma\left((T,v_{\point})\right)$ on $\varrho^{\lambda(v_{\point})}(\calP_n)$ is defined by
\begin{itemize}
    \item $\omega_\east = \{u_{j} : \lambda(x) = j \text{ for some internal vertex } x\}$, and
    \item $\omega_\north = \{v_{j} : \lambda(a) = j \text{ for some auxiliary vertex } a\}$.
\end{itemize}
\end{thm}

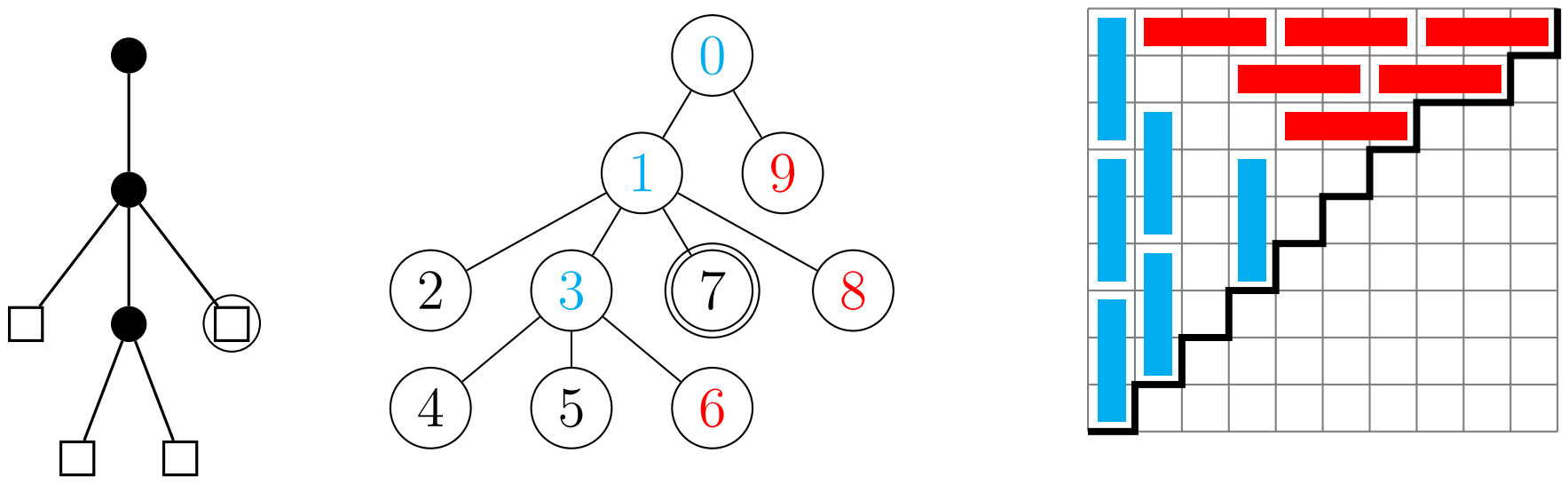
\begin{figure}[ht]
\centering
\begin{tikzpicture}[scale=.35,
   labelled/.style={
    circle, draw, minimum size=6mm, inner sep=0pt
  },
  dot/.style={
    circle, fill=black, inner sep=0pt, minimum size=7pt
  },
  level distance=25mm,
  sibling distance=30mm,
  edge from parent/.style={draw}
]
\begin{scope}[scale=4]
\begin{forest}
  baseline,
  tree,
  for tree={
    l=10mm,l sep=0pt,s sep=5mm,
    minimum size = 7pt,
  }
  [
    [ [] [ [] []  ] [] ]
  ]
\end{forest}
\draw (-.095,-1.425) circle[radius=0.15cm];
\end{scope}
        
    \begin{scope}[shift={(15,0)}, every node/.style={labelled}]
        \node {\textcolor{cyan}{0}}
          child {node {\textcolor{cyan}{1}}
            child {node {2}}
            child {node {\textcolor{cyan}{3}}
                child {node {4}}
                child {node {5}}
                child {node {\textcolor{red}{6}}}
            }
            child {node {7}}
            child {node {\textcolor{red}{8}}}
          }
          child {node {\textcolor{red}{9}}};
          \draw (0,-5) circle[radius=1cm];
    \end{scope}

    \begin{scope}[shift={(23,-8)}]
        \draw[step=1,color=gray] (0,0) grid (10,9);
        \draw[line width=1.5,color=black] (0,0)--(1,0)--(1,1)--(2,1)--(2,2)--(3,2)--(3,3)--(4,3)--(4,4)--(5,4)--(5,5)--(6,5)--(6,6)--(7,6)--(7,7)--(9,7)--(9,8)--(10,8)--(10,9);
        \brect{0}{0}
        \brect{0}{3}
        \brect{0}{6}
        \brect{1}{1}
        \brect{1}{4}
        \brect{3}{3}
        \rrect{10}{9}
        \rrect{7}{9}
        \rrect{4}{9}
        \rrect{9}{8}
        \rrect{6}{8}
        \rrect{7}{7}
    \end{scope}
    \end{tikzpicture} 
    \caption{For a vertex-pointed $1$-leafy tree $T$ (left), we depict the tree $\aux(T)$ (middle) where $v_{\point}$ is circled in each.  On the right, we have the corresponding tight grading $\omega$, where $\omega_\east = \{u_0,u_1,u_3\}$ and $\omega_\north = \{v_6,v_8,v_9\}$. Note that the tight grading is $7$-framed, where $7$ is the label of the distinguished vertex in $\aux(T)$.} 
    \label{fig: leafy tree to tight grading}

\end{figure}

\begin{rem}
    The inverse map to $\widehat \Gamma$ can be readily constructed as follows.  Looking each horizontal edge $u_i \in \varrho^f(\calP_n)$, create a word $w_1w_2\cdots w_n$ on the $3$-letter alphabet $\{(,),\bullet\}$ as follows:
    $$w_i =
    \begin{cases}
         (  & \text{ if } \omega(u_i) > 0,\\
         )  & \text{ if } \omega(v_i) > 0, \text{ where $v_i$ is a north step immediately following } u_i,\\
         \bullet & \text{ otherwise.}
    \end{cases}
    $$
    The leafy tree $T$ is the plane tree associated to $w_1w_2\cdots w_n$, and $v_{\point}$ is the vertex corresponding to the unique horizontal edge that is followed by another horizontal edge.
    
\end{rem}

\subsection{\texorpdfstring{Bijection between $0$-framed tight gradings and leafy trees}{Bijection between 0-framed tight gradings and leafy trees}}

In order to construct a bijection involving leafy trees rather than the vertex-pointed leafy trees, we refine the tight gradings by their frame.  Let $\TG{m}{n}{f}$ denote the subset of tight gradings in $\TGS{m}{n}$ that are $f$-framed.  

Note that the set of leafy trees in $\LTS{m}{n}$ is naturally identified with the subset of $\PLTS{m}{n}$ where the distinguished vertex is the root, in which case we say the tree is \emph{root-pointed}.  Equivalently, the root-pointed leafy trees in $\PLTS{m}{n}$ are those where $\lambda(v_{\point}) = 0$.  Thus under the map $\widehat{\Gamma}$, the subset of root-pointed leafy trees in $\PLTS{m}{n}$ corresponds precisely to $\TG{m}{n}{0}$. 

These identifications together yield a direct bijection from leafy trees to $0$-framed tight gradings.

\begin{cor}
    There is an explicit bijection 
    \begin{equation}\label{eq:leafy to tight restricted}
    \Gamma \colon \LTS{m}{n} \to \TG{m}{n}{0}\,    
    \end{equation}
    taking $T \in \LTS{m}{n}$ with root $v_0$ to $\widehat\Gamma\left((T,v_0)\right)$, as constructed in \Cref{thm: leafy to tight bij}. 
\end{cor}

\begin{exmp}
    Let $T$ be the $2$-leafy tree at the left of \Cref{fig: leafy tree to tight grading} (without a distinguished vertex).  Then the corresponding tight grading $\Gamma(T)$ is depicted on the left of \Cref{fig: tight grading examples}.
\end{exmp}

\section*{Acknowledgements}
The authors are grateful to Hyun Kyu Kim and Korea Institute for Advanced Study, where several key ideas were developed when AB and LM visited KL in June 2026. The authors also thank Fran\c{c}ois Bergeron, Mireille Bousquet-M\'elou, Fr\'ed\'eric Chapoton, \'Eric Fusy, Vincent Pilaud, Markus Reineke, and Hugh Thomas for helpful discussions and comments.  AB was partially supported by NSF grant DMS-2503411.

{\bf AI disclosure.}  The authors formulated the definition of leafy trees, inspired by their study of tight gradings, without AI assistance.  After this, the only notable mathematical AI contribution came from querying GPT-5.6 in July 2026 for ways to associate a Dyck path to a $1$-leafy tree.  Among the suggestions were the constructions we later generalized to define the bijection $\Phi$.  The proofs (to be presented in a later full version of this article) were discovered without AI assistance.  AI was used elsewhere in reviewing exposition, performing literature searches, doing numerical experiments, and generating figures.

\bibliographystyle{amsplain}
\bibliography{arxiv_announcement.bib} 
\addtocontents{toc}{\protect\setcounter{tocdepth}{1}}

\AtEndDocument{%
  \includepdf[
    pages=1,
    fitpaper=true,
    pagecommand={\thispagestyle{empty}}
  ]{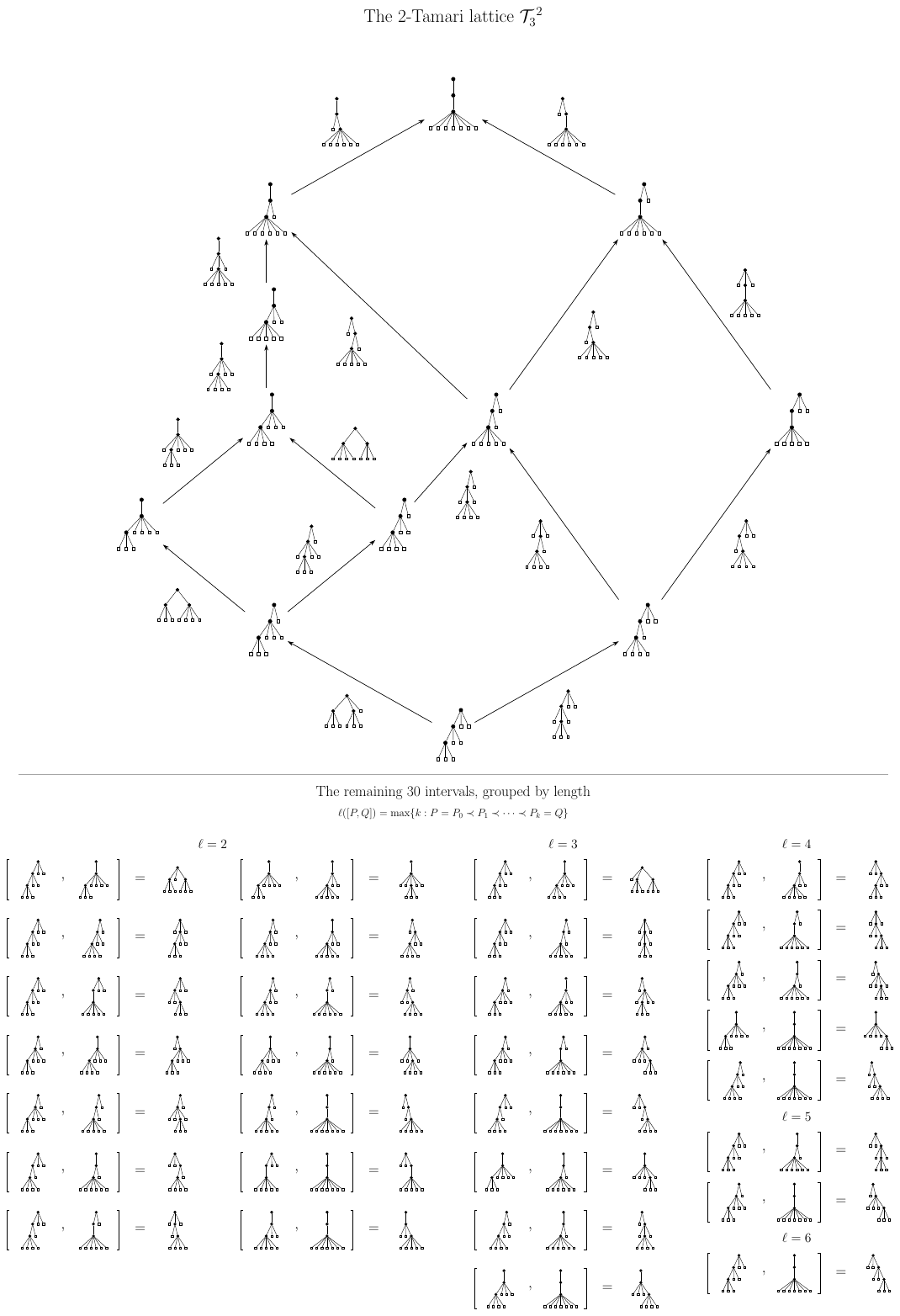}%
}

\end{document}